\documentclass{paper}

\usepackage{csquotes}

\newcommand{\E}{\mathcal{E}}

\renewcommand{\Db}{\D^\mathrm{b}} 
\renewcommand{\Pic}{\operatorname{Pic}} 

\newcommand{\cK}{\mathcal K} 
\newcommand{\cL}{\mathcal L} 

\newcommand{\cG}{\mathcal G}
\newcommand{\R}{\mathbf R} 

\newcommand{\Ample}{\operatorname{Ample}}

\newcommand{\PG}[1]{{\color{PineGreen}#1} }

\newtheorem{theorem*}{Theorem}

\title{Splitting of vector bundles on toric varieties}

\author{Mahrud Sayrafi}
\address{School of Mathematics, University of Minnesota, Minneapolis, Minnesota, 55455}
\email{\href{mailto:mahrud@umn.edu}{mahrud@umn.edu}}
\urladdr{\url{https://math.umn.edu/~mahrud/}}


\begin{document}

\begin{abstract}
  We prove a Horrocks-type splitting criterion for arbitrary smooth projective toric varieties under an additional hypothesis similar to the case of products of projective spaces by Eisenbud--Erman--Schreyer.
\end{abstract}

\maketitle

\section{Introduction}

The study of algebraic vector bundles as a rich source of high dimensional varieties, nonsingular subvarieties, and moduli problems is classical in algebraic geometry \cite{Hartshorne74,Hartshorne79}. Moreover, the equivalence of the categories of algebraic and holomorphic vector bundles on a complex algebraic variety connects this study to mathematical physics.

One central problem here is determining the indecomposability of a given vector bundle. \linebreak By a famous result of Horrocks, a vector bundle $\E$ on $\PPn$ splits as a direct sum of line bundles if and only if $H^q(\PPn,\E(a)) = 0$ for all $q = 1,\dots,n-1$ and all twists $\OO(a)\in\Pic\PPn$ \cite{Horrocks64}. Barth and Hulek gave an inductive proof of this by restricting to a linear subspace $\PP{n-1}$ and using Grothendieck's theorem for $\PP1$ as the base \cite[Lem.~1]{BH78}.

Horrocks' splitting criterion inspired similar criteria for splitting of vector bundles over other classes of varieties: products of projective spaces \cite{CM05,EES15}, Grassmannians and quadrics \cite{Ott89}, rank 2 vector bundles on Hirzebruch surfaces \cite{Buchdahl87,AM11,FM11,Yas15}, and Segre--Veronese varieties \cite{Schreyer22}, among others. See \cite{Ottaviani24} for a recent survey.

We prove an analogous splitting criterion for vector bundles on smooth projective toric varieties, under an additional hypothesis similar to that of Eisenbud--Erman--Schreyer's criterion for products of projective spaces \cite[Thm.~7.2]{EES15} and the recent joint work with Brown for the Picard rank 2 case in \cite[Thm.~1.5]{BS24}. 

\begin{theorem*}\label{main-thm-horrocks}
  Suppose $\E$ is a vector bundle on a smooth projective toric variety $X$ and $\E' = \oplus_{i=1}^n \OO(D_i)^{r_i}$ is a sum of line bundles on $X$ such that $D_{i+1}-D_i$ is ample for $0<i<n$. \linebreak
  If $H^q(X,\E\otimes\cL) = H^q(X,\E'\otimes\cL)$ for all $q\geq0$ and $\cL\in\Pic X$, then $\E\cong\E'$.
\end{theorem*}

The main ingredient in the proof of Theorem~\ref{main-thm-horrocks} is a construction of resolutions of toric subvarieties by line bundles due to Hanlon--Hicks--Lazarev \cite{HHL24} (c.f.~\cite{FH23,Anderson23,BE24a}). The proof consists of a Beilinson-type spectral sequence which computes the Fourier--Mukai transform corresponding to a resolution of the diagonal. Similar ideas have been used to great success in \cite{CM05,FM11,AM11,EES15,BS24}.

In our case, a significant obstacle is introduced by the difference between the nef and effective cones for arbitrary toric varieties. In all previous incarnations of the criterion, either the nef and effective cones are identical or the Picard rank is low (one or two). Without either of these restrictions, the analysis of the cohomology of line bundles requires new ideas.

\subsection*{Outline}

We begin in \S\ref{sec:splitting-recipe} with a recipe for proving Horrocks-type splitting criteria for arbitrary smooth projective varieties, which illustrates the proof. Then in \S\ref{sec:toric-horrocks} we prove Theorem~\ref{main-thm-horrocks}.

\subsection*{Acknowledgments}

The author thanks Christine Berkesch, Lauren Cranton Heller, Jay Yang, Gregory Smith, Michael Brown, Daniel Erman, Andrew Hanlon, Jeff Hicks, Oleg Lazarev, David Favero, Jesse Huang, Devlin Mallory, and David Eisenbud for valuable conversations. This paper was written while the author was supported by the Doctoral Dissertation Fellowship at the University of Minnesota.

\section{A general recipe for splitting criteria}\label{sec:splitting-recipe}

Let $X$ be a smooth projective variety with a resolution of the diagonal $\cK$ and $\E$ a coherent sheaf on $X$. Similar to the case in \cite[\S4]{BS24}, we use a Fourier--Mukai functor with kernel $\cK$ to construct a monad which is quasi-isomorphic to $\E$ and whose terms are prescribed by the terms of $\cK$ with ranks given by sheaf cohomology of twists of $\E$. The recipe for the splitting criteria proved in this section is a consequence of appropriate vanishing of the terms of this spectral sequence.

The diagonal embedding $X \to X \times X$ defines a closed subscheme $\Delta \subset X \times X$. Let $\pi_1$ and $\pi_2$ denote the natural projections of $X \times X$ onto $X$ and for the rest of this paper suppose $\cK$ is a locally free resolution for $\OO_\Delta$, the structure sheaf of $\Delta$, with terms given as direct sums of sheaves of the form $\cG\boxtimes\cL\coloneq \pi_1^*\cG\otimes\pi_2^*\cL$, where $\cG$ is a locally free sheaf and $\cL=\OO(E)$ is a line bundle corresponding to a divisor $E$ on $X$.

The Fourier--Mukai transform with kernel $\cK$ is the composition of functors:
\[\Phi_\cK\colon \Db(X) \xto{\pi_1^*} \Db(X \times X) \xto{\cdot\;\otimes\cK} \Db(X \times X) \xto{\R\pi_{2*}} \Db(X). \]
In particular, $\Phi_\cK$ is the identity functor on the derived categories, meaning that $\Phi_\cK(\E)$ will be quasi-isomorphic to $\E$. We evaluate the last functor, derived pushforward, by resolving the first term of each box product with a \v{C}ech complex to obtain a spectral sequence
\[ E_1^{-p,q} = \R^q\pi_{2*}(\pi_1^*\E \otimes \cK_p) = \bigoplus_i \cG_i \otimes H^q(X, \cL_i) \Rightarrow \R^{p-q}\pi_{2*}(\pi_1^*\E \otimes \cK) \cong \begin{cases} \E & i = j \\ 0 & i \neq j, \end{cases} \]
where the direct sum ranges over summands $\cG_i\boxtimes\cL_i$ of $\cK_p$ (c.f.~\cite[\S3.3]{BCHS21} and \cite[\S4]{BS24}).

\begin{figure}[h]\label{fig:E1}
  \begin{equation}\label{eq:E1}
    \begin{tikzcd}[column sep=small, row sep=tiny]
      {} & \vdots & \vdots & \vdots & {} \\[-5pt]
      {} & {\R^2{\pi_2}_*(\pi_1^*\E\otimes\cK_0)} & {\R^2{\pi_2}_*(\pi_1^*\E\otimes\cK_1)} &     {\R^2{\pi_2}_*(\pi_1^*\E\otimes\cK_2)} & \cdots \\
      {} & {\R^1{\pi_2}_*(\pi_1^*\E\otimes\cK_0)} & {\R^1{\pi_2}_*(\pi_1^*\E\otimes\cK_1)} & \red{\R^1{\pi_2}_*(\pi_1^*\E\otimes\cK_2)} & \cdots \\
      {} &  \PG{{\pi_2}_*(\pi_1^*\E\otimes\cK_0)} & \red{{\pi_2}_*(\pi_1^*\E\otimes\cK_1)} &     {    {\pi_2}_*(\pi_1^*\E\otimes\cK_2)} & \cdots \\
      {} & {} & {}
      \arrow[from=4-4, to=4-3]
      \arrow[from=4-3, to=4-2]
      \arrow[from=4-5, to=4-4]
      \arrow[from=3-4, to=3-3]
      \arrow[from=3-5, to=3-4]
      \arrow[from=2-5, to=2-4]
      \arrow[from=2-4, to=2-3]
      \arrow[from=2-3, to=2-2]
      \arrow[from=3-3, to=3-2]
      \arrow["{t}"{description, sloped, pos=1.3}, crossing over, shift left=17, shorten <=-28pt, shorten >=-18pt, dashed, no head, from=4-2, to=1-2]
      \arrow["{s}"{description, sloped, pos=1.1}, shift right=5, shorten <=-87pt, shorten >=-20pt, dashed, no head, from=4-2, to=4-5]
      \arrow["{k=2}"{description, sloped, pos=0.70}, shift left=4, shorten >=70pt, dotted, no head, from=3-5, to=5-3]
      \arrow["{k=1}"{description, sloped, pos=0.81}, shift left=4, shorten >=70pt, dotted, no head, from=2-5, to=5-2]
      \arrow["{k=0}"{description, sloped, pos=0.95}, shift left=1, shorten >=30pt, dotted, no head, from=1-5, to=5-1]
      \arrow[shift right=2, shorten <=25pt, shorten >=23pt, dotted, no head, from=1-3, to=3-1]
      \arrow[shift right=2, shorten <=25pt, shorten >=23pt, dotted, no head, from=1-4, to=4-1]
    \end{tikzcd}
  \end{equation}
\end{figure}

\begin{definition}\label{def:support-cone}
  For a convex cone $\mathcal A \subset \Pic X$, say $\cK$ is \emph{cohomologically supported in $\mathcal A$} if any summand $\cG\boxtimes\OO(E)$ of $\cK_p$ has $H^q(X, \OO(E-D)) = 0$ for all $q<p$ and $D\in\mathcal A$. \linebreak This means the terms along dotted diagonals with $k>0$ in \eqref{eq:E1} vanish when $\E = \OO(-D)$.
\end{definition}

For instance, Beilinson's resolution of the diagonal for $\PPn$ \cite{Beilinson78}, its variant for products of projective spaces (e.g.~\cite{BES20}), and the resolutions constructed in \cite{BS24} are all cohomologically supported in $\Nef X\subset\Pic X$. Note that only the Picard group for $\PPn$ has a total ordering, hence any direct sum of line bundles on $\PPn$ can be increasingly ordered by $\Nef\PPn$.


\begin{lemma}\label{lem:splitting}
  Let $\E$ be a coherent sheaf on a smooth projective variety $X$ with a resolution of the diagonal $\cK$ such that $\Phi_\cK(\OO_X) = \OO_X$. Consider the spectral sequence $E_1^{-p,q} \Rightarrow \E$ above.
  \begin{enumerate}[ref={\thelemma(\alph*)}]
  \item\label{lem:reds-vanish} If $E_1^{-p-1,p} = 0$ for all $p$ (the red terms vanish), then $E_1^{0,0}$ is a direct summand of $\E$.
  \item\label{lem:cone-vanish} If $\cK$ is supported in $\mathcal A$ and $\E = \bigoplus\OO(D_i)$ with $-D_i\in\mathcal A$, then $E_1^{-p-1,p} = 0$ for all $p$.
  \end{enumerate}
\end{lemma}
\begin{proof}
  The proof of the first part is identical to \cite[Lem.~4.1]{BS24} and \cite[Lem.~7.3]{EES15}.\linebreak Using \cite[Lem.~3.5]{EFS03}, there exists a complex with terms that are the same as the totalization $\Tot(E_1)$ (along the dotted diagonals in \eqref{eq:E1}) which is quasi-isomorphic to $\E$. The vanishing of the first term of the totalization (colored in red) implies that all differentials with source or target $E_r^{0,0}$ are zero, therefore $E_1^{0,0}$ is a summand of $E_\infty^{0,0} = \E$.

  The second part immediately follows from Definition~\ref{def:support-cone}, as $\Phi_\cK$ commutes with direct sums and the totalization $\Tot(E_1)$ corresponding to $\Phi_\cK(\OO(-D))$ for any divisor $D\in\mathcal A$ is zero in positive homological degrees (i.e., terms along dotted diagonals with $k>0$ vanish).
\end{proof}

\begin{remark}
  The additional hypothesis present in the splitting criteria proved in \cite{EES15,BS24} avoids the problem of missing a total ordering in the case of higher Picard rank by proving a criterion for a smaller cone $\mathcal A = \Nef X$. Nevertheless, to date we do not know whether such hypotheses are necessary. In contrast, for varieties such as Hirzebruch surfaces, there are resolutions of the diagonal supported in $\Eff X$, yielding a splitting criterion which requires a strictly weaker additional hypothesis.
\end{remark}


\begin{proposition}\label{prop:splitting}
  Suppose $X$ is a smooth projective variety $X$ with a locally free resolution of the diagonal $\cK$ such that $\cK$ is cohomologically supported in $\mathcal A$ and $\Phi_\cK(\OO_X) = \OO_X$. Let $\E$ be a vector bundle and $\E' = \oplus_{i=1}^n \OO(D_i)^{r_i}$ a sum of line bundles on $X$ such that $D_{i+1}-D_i\in\mathcal A$ for $0<i<n$.
  If $H^q(X, \E\otimes\cL) = H^q(X, \E'\otimes\cL)$ for all $q\geq0$ and $\cL\in\Pic X$, then $\E\cong\E'$.
\end{proposition}
\begin{proof}
  Following a similar road map as \cite[Thm.~1.5]{BS24} and \cite[Thm.~7.2]{EES15}, twist $\E$ and $\E'$ by the highest line bundle $\OO(-D_n)$ so that without loss of generality we can assume $\E' = \oplus_{i=1}^{n-1} \OO(D_i)^{r_i} \oplus \OO_X^{r_n}$. Let $E_1(\E)$ denote the spectral sequence corresponding to $\Phi_\cK(\E)$.

  By hypothesis, $E_1(\E)$ and $E_1(\E')$ have the same terms, so $E_1^{0,0}(\E) = E_1^{0,0}(\E') = \OO_X^{r_n}$. Since $-D_i\in\mathcal A$ for all $i$, Lemma~\ref{lem:cone-vanish} implies that $E_1^{-p-1,p}(\E) = E_1^{-p-1,p}(\E') = 0$ for all $p$. Using Lemma~\ref{lem:reds-vanish}, the term $E_1^{0,0}(\E) = \OO_X^{r_n}$ is a summand of $\E$. Induction on the complement of $\OO_X^{r_n}$ in $\E$ and $\E'$ finishes the proof.
\end{proof}

\begin{remark}
  Proposition~\ref{prop:splitting}, combined with resolutions of the diagonal constructed by Beilinson \cite{Beilinson78} and Kapranov \cite{Kapranov88}, recovers the splitting criteria for $\PPn$ by Horrocks and Grassmannians and quadrics by Ottaviani \cite{Ott89}, respectively. Further, one can show that given resolutions of the diagonal $\cK$ and $\cK'$ supported in $\mathcal A$ and $\mathcal A'$ for $X$ and $X'$,\linebreak respectively, $\cK\boxtimes\cK'$ is a resolution of the diagonal for $X \times X'$ supported in $\mathcal A \times \mathcal A'$.\linebreak In particular, the splitting criterion \cite[Thm.~7.2]{EES15} for products of projective spaces can be recovered from Beilinson's resolution for $\PPn$, similarly for products of Grassmannians, etc.
\end{remark}


\pagebreak
\section{Smooth Projective Toric Varieties}\label{sec:toric-horrocks}

In \cite{HHL24}, Hanlon, Hicks, and Lazarev construct resolutions of toric subvarieties by line bundles on a smooth toric variety $X$. The case of interest here is the diagonal subvariety, where the resolution of the structure sheaf of the diagonal $\cK$ consists of line bundles from the Thomsen collection on $X \times X$ (compare with \cite{BE24a}).

We will need the following technical lemma on properties of the terms of $\cK$. We denote by $N_X$ and its dual $M_X = \Hom(N_X, \ZZ)$ the lattices of one-parameter subgroups and characters of the torus, respectively, and $N_{X,\RR}$ or $M_{X,\RR}$ the corresponding real vector spaces. Further, for a ray $\rho\in\Sigma(1)$ denote by $D_\rho$ the corresponding torus-invariant prime divisor on $X$.

\begin{lemma}\label{lem:terms}
  Suppose $\OO(E')\boxtimes\OO(E)$ is a summand of $\cK_p$ constructed as in \cite{HHL24}.
  \begin{enumerate}[label=\normalfont{(\Roman*).}, ref=\Roman*]
  \item\label{item:fact-1} The divisor $-E$ is an effective Cartier divisor on $X$; that is,
    \[ E = -\sum d_\rho D_\rho \quad \text{ for } d_\rho\in\ZZ_{\geq0} \text{ and } \rho\in\Sigma(1). \]
  \item\label{item:fact-2} The bundle $\OO(E)$ is a summand of a high toric Frobenius pushforward of $\OO_X$; that is, there is a Cartier $\QQ$-divisor $\tilde E$ linearly equivalent to $E$ such that
    \[ \tilde E = -\sum c_\rho D_\rho \quad \text{ for } c_\rho\in[0,1) \text{ and } \rho\in\Sigma(1). \]
  \item\label{item:fact-3} The dimension of the polyhedron $P_{-E}$ of the divisor $-E$ is at least $p$; that is,
    \[ \text{if } \OO(E')\boxtimes\OO(E) \text{ is a summand of } \cK_p, \text{ then } p\leq\dim P_{-E}. \]
  \end{enumerate}
\end{lemma}
\begin{proof}
  The Thomsen collection for $X \times X$ consists of products of bundles from the Thomsen collection for $X$ \cite[Rem.~1.3]{HHL24}, hence the first two properties follow from $\OO(E)$ being in the Thomsen collection for $X$ \cite[\S5]{HHL24}.

  The third point is more subtle, as it implies that not all line bundles from the Thomsen collection for $X \times X$ appear in $\cK$, and only a subset may appear in any given term $\cK_p$.\linebreak The diagonal embedding is induced by an inclusion of lattices $\bar\phi\colon N_X \to N_{X \times X}$. The dual map on the character lattices $\bar\phi^*\colon M_{X \times X} \to M_X$ induces a short exact sequence of real tori:
  \[ 0 \to L_\RR \to M_{X \times X, \RR}/M_{X \times X} \to M_{X, \RR}/M_X \to 0. \]
  The construction of $\cK$ begins with a stratification $S$ of the ambient real torus $M_{X \times X, \RR}/M_{X \times X}$ labeled by divisors on $X \times X$ introduced by Bondal \cite{Bondal06} (c.f.~\cite[\S3.4]{HHL24}). In the case of the diagonal, the dimension of the kernel $L_\RR \cong M_{X,\RR}/M_X$ equals $\dim X$ and the kernel inherits the stratification $S$ (c.f.~\cite[Exa.~3.13]{HHL24}).

  We need to show that if $\OO(E')\boxtimes\OO(E)$ is a summand of $\cK_p$ then $p\leq\dim P_{-E}$. It follows from the construction of $\cK$ in \cite[eq.~(18)]{HHL24} that the line bundle summands of $\cK_p$ correspond to labels on the $p$-dimensional strata in $S$. In particular, the stratification on $L_\RR$ is the same as the stratification on $M_{X,\RR}$ when constructing the resolution of a point on $X$, only in that case the labels are divisors on $X$ rather than $X \times X$. Specifically, if a $p$-dimensional stratum $S_\sigma$ has label $\OO(E')\boxtimes\OO(E)$ in $\cK_p$, then $S_\sigma$ has label $\OO(E)$ in the resolution of a point on $X$. We will use this correspondence to bound the dimension of $P_{-E}$.

  Given a line bundle $\OO(E)$ in the Thomsen collection, let $S_{-E}$ denote the union of strata with that label in Bondal's stratification on $M_{X,\RR}/M_X$. It follows from \cite[Lem.~5.6]{FH22} that $S_{-E} = (P_{-E} \setminus \bigcup_{\rho\in\Sigma(1)} P_{-E-D_\rho}) / M_X$. Since $-E$ is effective, $P_{-E}$ is nonempty, and because the section polytopes are closed, $S_{-E}$ is open and $\dim S_{-E} = \dim P_{-E}$. Putting this all together, any $p$-dimensional stratum $S_\sigma$ which is labeled by the line bundle $\OO(E')\boxtimes\OO(E)$ in $\cK_p$ must satisfy $p = \dim S_\sigma \leq \dim S_{-E} = \dim P_{-E}$.
\end{proof}

\begin{remark}
  While the stratifications considered in \cite[\S3.4]{HHL24} and \cite[\S5]{FH22} are both versions of the stratification studied by Bondal in \cite{Bondal06}, they have a subtle difference: the union of the strata with the same label in \cite{HHL24} is the unique strata with that label in \cite{FH22}, which is contractible by \cite[Lem.~5.6]{FH22}.
\end{remark}


In order to use Lemma~\ref{lem:splitting}, we need the following analysis of the support of $\cK$.

\begin{proposition}\label{prop:support}
  The resolution of the diagonal $\cK$ is cohomologically supported in $\Ample(X)$.
\end{proposition}
\begin{proof}
  \newcommand{\eps}{\epsilon}
  Suppose $\OO(E')\boxtimes\OO(E)$ is a summand of $\cK_p$. We show that:
  \[ H^q(X, \OO(E-D)) = 0 \text{ for } q<p \text{ and any ample divisor } D. \]
  First, using $d_\rho$ and $c_\rho$ from \eqref{item:fact-1} and \eqref{item:fact-2}, since $\lceil d_\rho-(1-\eps)c_\rho\rceil = d_\rho$ for any $\eps\in[0,1]$,
  \begin{align*}
    -\lceil D + (1-\eps)\tilde E - E\rceil
    &= -\lceil D + \sum(d_\rho - (1-\eps)c_\rho) D_\rho\rceil \\
    &= -\lceil D + \sum d_\rho D_\rho\rceil = -\lceil D - E\rceil = E - D.
  \intertext{Second, by linear equivalence in \eqref{item:fact-2} it follows that $\tilde E - E \sim 0$. Hence,}
    D + (1-\eps)\tilde E - E
    &= D + (\tilde E - E) - \eps\tilde E \\
    &\sim D - \eps\tilde E,
  \intertext{which, for sufficiently small $\eps$, is ample, and hence nef, because $D$ is ample by hypothesis.
      Third, since both $D$ and $-\eps\tilde E$ are effective,}
  \dim P_{D + (1-\eps)\tilde E - E} &= \dim P_{D -\eps\tilde E} \\
  &\geq \dim P_{-\eps\tilde E} = \dim P_{-\tilde E} = \dim P_{-E}.
  \end{align*}
  Hence by Batyrev--Borisov Vanishing (see \cite[Thm.~9.3.5(b)]{CLS2011}), 
  \[ H^q(X, \OO(E-D)) = H^q(X, \OO(-\lceil D + (1-\eps)\tilde E - E\rceil)) = 0 \text{ for all } q<\dim P_{-E}. \]
  Therefore by \eqref{item:fact-3}, the weaker vanishing $H^q(X, \OO(E-D)) = 0$ for $q<p$ holds.
\end{proof}

The proof of the main theorem is a direct application of the recipe from Section~\ref{sec:splitting-recipe}.

\begin{proof}[Proof of Theorem~\ref{main-thm-horrocks}]
  Since $\cK$ is cohomologically supported in $\Ample(X)$ by Proposition~\ref{prop:support} and the consecutive differences $D_{i+1}-D_i$ are ample by hypothesis, the proof follows immediately using the recipe in Proposition~\ref{prop:splitting}.
\end{proof}


\bibliography{references}

\end{document}